%


\ifx\UsualIsLoaded\undefined
\let\UsualIsLoaded=\relax		

\input psfig	

\font\fourteenrm=cmr12  scaled \magstep1
\font\fourteenbf=cmbx12 scaled \magstep1
\font\fourteentt=cmtt12 scaled \magstep1
\font\fourteensl=cmsl12 scaled \magstep1
\font\fourteensy=cmsy10 scaled \magstep2
\font\fourteeni=cmmi12  scaled \magstep1
\font\fourteenit=cmti12 scaled \magstep1
\font\fourteensc=cmcsc10 scaled \magstep2
\font\fourteenbit=cmssi12 scaled \magstep1
\font\fourteenBbb=msym10 scaled \magstep2

\font\twelverm=cmr12
\font\twelvebf=cmbx12
\font\twelvett=cmtt12
\font\twelvesl=cmsl12
\font\twelvesy=cmsy10 scaled \magstep1
\font\twelvei=cmmi12
\font\twelveit=cmti12
\font\twelvesc=cmcsc10 scaled \magstep1
\font\twelvebit=cmssi12
\font\twelveBbb=msym10 scaled \magstep1
\newfam\Bbbfam
\font\tenrm=cmr10
\font\tenbf=cmb10 
\font\tentt=cmtt10 
\font\tensl=cmsl10 
\font\tensy=cmsy10 
\font\teni=cmmi10 
\font\tenit=cmti10 
\font\tensc=cmcsc10 
\font\tenbit=cmssi10
\font\tenBbb=msym10
\font\ninei=cmmi9 
\font\ninerm=cmr9

\font\ninesy=cmsy9 
\font\nineBbb=msym9
\font\eighti=cmmi8
\font\eightrm=cmr8

\font\eightsy=cmsy8
\font\eightBbb=msym8
\font\seveni=cmmi7 
\font\sevenrm=cmr7
\font\sevensy=cmsy7
\font\sevenBbb=msym7
 
%
%

%

\def\tenpoint{
\def\rm{\fam0\tenrm}
\textfont0=\tenrm \scriptfont0=\eightrm \scriptscriptfont0=\sevenrm
\textfont1=\teni \scriptfont1=\eighti \scriptscriptfont1=\seveni
\textfont2=\tensy \scriptfont2=\eightsy \scriptscriptfont2=\sevensy
\textfont3=\tenex \scriptfont3=\tenex \scriptscriptfont3=\tenex

\textfont\itfam=\tenit
\def\it{\fam\itfam\tenit}

\textfont\slfam=\tensl
\def\sl{\fam\slfam\tensl}

\textfont\bffam=\tenbf
\def\bf{\fam\bffam\tenbf}

\textfont\ttfam=\tentt
\def\tt{\fam\ttfam\tentt}

\def\sc{\tensc}

\def\bit{\tenbit}

\def\Bbb{\fam\Bbbfam\twelveBbb}
\textfont\Bbbfam=\tenBbb
\scriptfont\Bbbfam=\eightBbb 
\scriptscriptfont\Bbbfam=\sevenBbb
 
\normalbaselineskip=12pt

\setbox\strutbox=\hbox{\vrule height10pt depth4pt width0pt}%
\normalbaselines\rm}


\def\twelvepoint{
\def\rm{\fam0\twelverm}
\textfont0=\twelverm \scriptfont0=\ninerm \scriptscriptfont0=\sevenrm
\textfont1=\twelvei \scriptfont1=\ninei \scriptscriptfont1=\seveni
\textfont2=\twelvesy \scriptfont2=\ninesy \scriptscriptfont2=\sevensy
\textfont3=\tenex \scriptfont3=\tenex \scriptscriptfont3=\tenex

\textfont\itfam=\twelveit
\def\it{\fam\itfam\twelveit}

\textfont\slfam=\twelvesl
\def\sl{\fam\slfam\twelvesl}

\textfont\bffam=\twelvebf
\def\bf{\fam\bffam\twelvebf}

\textfont\ttfam=\twelvett
\def\tt{\fam\ttfam\twelvett}

\def\sc{\twelvesc}
\def\bit{\twelvebit}

\def\Bbb{\fam\Bbbfam\twelveBbb}
\textfont\Bbbfam=\twelveBbb
\scriptfont\Bbbfam=\nineBbb
\scriptscriptfont\Bbbfam=\sevenBbb
 
\normalbaselineskip=14pt

\setbox\strutbox=\hbox{\vrule height10pt depth4pt width0pt}%
\normalbaselines\rm}


\def\fourteenpoint{
\def\rm{\fam0\fourteenrm}
\textfont0=\fourteenrm \scriptfont0=\twelverm \scriptscriptfont0=\tenrm
\textfont1=\fourteeni \scriptfont1=\twelvei \scriptscriptfont1=\teni
\textfont2=\fourteensy \scriptfont2=\twelvesy \scriptscriptfont2=\tensy
\textfont3=\tenex \scriptfont3=\tenex \scriptscriptfont3=\tenex

\textfont\itfam=\fourteenit
\def\it{\fam\itfam\fourteenit}

\textfont\slfam=\fourteensl
\def\sl{\fam\slfam\fourteensl}

\textfont\bffam=\fourteenbf
\def\bf{\fam\bffam\fourteenbf}

\textfont\ttfam=\fourteentt
\def\tt{\fam\ttfam\fourteentt}

\def\sc{\fourteensc}
\def\bit{\fourteenbit}

\def\Bbb{\fam\Bbbfam\twelveBbb}
\textfont\Bbbfam=\fourteenBbb
\scriptfont\Bbbfam=\tenBbb
\scriptscriptfont\Bbbfam=\eightBbb

\normalbaselineskip=16pt

\setbox\strutbox=\hbox{\vrule height10pt depth4pt width0pt}%
\normalbaselines\rm}

%
\twelvepoint

\abovedisplayskip 14pt plus 3pt minus 10pt%
\belowdisplayskip 14pt plus 3pt minus 10pt%
\abovedisplayshortskip 0pt plus 3pt%
\belowdisplayshortskip 8pt plus 3pt minus 5pt%
\parskip 3pt plus 1.5pt
\hsize=6.5in
\vsize=8.9in

\fi			

\def\R{{\Bbb R}}

\def\N{{\Bbb N}}
\def\QED{\hfill\vbox{\hrule height1pt%
     \hbox{\vrule width1pt \kern4pt\vbox{\kern4pt\kern4pt}\kern4pt
	   \vrule width1pt} \hrule height1pt}\par\bigskip}

\long\gdef\ignore#1{}
\def\emptyset{/\kern -0.67em{\scriptstyle\bigcirc}}

\chardef \other = 12

\def\makeactive#1{\catcode`#1 = \active\ignorespaces}
{
	\makeactive\^^M %
	\gdef\obeywhitespace{%
		\makeactive\^^M %
		\let^^M = \newline %
		\aftergroup\removebox %
		\obeyspaces %
	}%
}
\def\newline{\par\indent}
\def\removebox{\setbox0=\lastbox}

\def\{}

\def\BoxIt#1#2{
	\vbox{\hrule
	\hbox{\vrule\kern#2\vbox{\kern#2#1\kern#2}\kern#2\vrule}
    \hrule}}

\def\R{{I\!\!R}}    \def\N{{I\!\!N}}  \def\e{\epsilon}

\centerline {\bf Cantor sets in the line:}
\centerline {\bf  scaling function and the
smoothness of the shiftmap}
\centerline {\it by F.Przytycki and F. Tangerman}
{\bf Abstract:} {\it Consider $d$ disjoint closed
subintervals of the unit
interval and consider an orientation preserving expanding map which
maps each of these subintervals
to the whole unit interval. The set of points where all iterates of
this
expanding map
are defined
is a Cantor set. Associated to the construction of this Cantor set is
the scaling
function which records the infinitely deep geometry of this Cantor
set.
This scaling function is an invariant of $C^1$ conjugation. We solve
the inverse problem posed by Dennis Sullivan:
given a scaling function, determine the maximal possible smoothness
of any expanding map which produces it.}\vskip .3 in
Consider the space $\Sigma_d\,=\,\lbrace\,1,..,d\,\rbrace^\N$, with its
standard shiftmap $\sigma$
$$\sigma(\alpha_1\alpha_2..)\,=\,(\alpha_2...)$$
Denote by $\sigma_i^{-1}$ the $d$ right-inverse of $\sigma$:
$$\sigma_i^{-1}(\alpha_1\alpha_2..)\,=\,(i\alpha_1\alpha_2..)$$
Our convention will be not to use separating comma's
in
strings of symbols.

$\Sigma_d$ with the product topology is a Cantor set.
Consider an embedding $h$ of the space
$\Sigma_d\,=\,\lbrace\,1,..,d\,\rbrace^\N$
into $\R$ with the standard order:
$$h(\alpha)\,>\,h(\beta)\;\;{\rm iff}\;\;\alpha_m\,>\,\beta_m$$
where $m$ is the first integer for which $\alpha_m\,\neq\,\beta_m$.
The image of $h$ is also a Cantor set. Denote by $f$ the induced
shiftmap on the image of $h$ and by $f_i^{-1}$ the $d$ right-inverses
of
$f$.
Let $r\,>\,1$.
We say that $h$ is $C^r$ if each of the right-inverses $f^{-1}_i$
have $C^r$ extensions to $\R$ which are contractions.
We say then that the Cantor set is $C^r$.\par
Every $C^{1+\epsilon}$ Cantor set has a scaling function, defined
below and there is a simple characterization of those functions
which are scaling functions for some $C^{1+\epsilon}$ Cantor set.
In this paper we describe those scaling functions which actually
have to $C^{k+\epsilon}$ realizations. Here $k$ is any integer
greater or equal to $1$ and $0\,<\,\epsilon\,\leq\,1$. We
follow the convention that $\epsilon\,=\,1$ means a Lipschitz
condition.

The theory for $r\,=\,1\,+\,\epsilon$ is essentially due to
Feigenbaum and Sullivan who introduced the scaling function.
It is defined in the following manner. Given an embedding $h$,
then the shiftmap allows a canonical definition of the image of $h$
as an intersection of nested collections of intervals. More
precisely, define for any finite sequence $(j_1..j_n)$
$I_{j_1..j_n}$ as the convex hull of
$h(\lbrace\,\alpha:\,\alpha_1\,=\,j_n,..,\alpha_n\,=\,j_1\,\rbrace)$
Note the order in which the indices occur.
Then for any $j_0$, $I_{j_0j_1..j_n}\,\subset\,I_{j_1..j_n}$
and the shiftmap maps $I_{j_1..j_n}$ to
$I_{j_1..j_{n-1}}$.
For the empty string, $I$ denotes the image of $h$.
The sets thus constructed are not intervals, but actually small
pieces of the image of $h$. It is however convenient to think of
them as intervals.\par
For any subset $J$ in the reals denote by $<J>$ its convex hull
and by
$|J|$ the length of its convex hull. We will in the remainder always
assume that $<I>$ is the unit interval [0,1].\par
Denote the set of finite strings $j_1,..,j_n$ of
length $n$ by $\Sigma_{d,n}^{dual}$.
The scaling function (ratio geometry) at level $n$ is a function
$S^n$:
$$S:\,\Sigma^{dual}_{d,n}\,\rightarrow\,(0,1)^{2d-1}$$
defined in the following manner.
For each $j_1..j_n$ $S(j_1..j_n)$ records the geometrical location
of the $d$ intervals
$\lbrace\,I_{j_0j_1..j_n}\,\rbrace_{j_0=1..d}$ in $I_{j_1..j_n}$
by the ratio's of lengths
of these $d$ intervals (first $d$ coordinates) and $d\,-\,1$ gaps (last
$d\,-\,1$ coordinates) to the length of $I_{j_1,..,j_n}$.
In particular for $j_0=1,..,d$ the $j_0-th$ coordinate of $S$ is
given by the following formula:
$$S(j_1..j_n)_{j_0}\,=\,{{|I_{j_0..j_n}|}\over{|I_{j_1..j_n}|}}$$
The sum of all ratio's of lengths equals one. Therefore $S$ actually
takes values in the $2d-2$ dimensional simplex $Simp_{2d-2}$ of
$(0,1)^{2d-1}$
where the sum of the coordinates equals 1.
Moreover lengths of intervals are determined by the scaling functions
at all levels:
$$|I_{j_1..j_n}|\,=\,\Pi_k\,S(j_{k+1}..j_n)_{j_k}\;\;(1)$$
Consider two finite sequences $j\,=\,j_1..j_n$ and
$j'\,=\,j'_1..j'_m$.
There is a canonical identification
between $I_j$
and $I_{j'}$ defined as follows. Let $j\,\cap\,j'$ be the longest
string which
agrees with both the beginning of $j$ and the beginning of $j'$.
Then suitable iterates of the shiftmap
map $I_j$ to $I_{j\,\cap\,j'}$ respectively $I_{j'}$ to
$I_{j\,\cap\,j'}$.(see
diagram)
$$\matrix{&I_{j'\cap j}&\cr
	  \nearrow& &\nwarrow\cr
	  I_{j'}& &I_j\cr}$$
The fundamental observation is
that if the embedding is
$C^{1+\epsilon}$ then the identification map is close to being
linear in the following precise sense. Define the nonlinearity
of a diffeomorphism $f$ on an interval as
$$\log\,\sup_{x,y,x\,\neq\,y}\,{{Df(x)}\over{Df(y)}}$$
Then the nonlinearity of the identification map can be
estimated from above in terms of the length of the intermediary
interval $I_{j\,\cap\,j'}$. But then
if $j\,\cap\,j'$ is long (i.e. $|I_{j\,\cap\,j'}|$ small),
the subdivision of $I_j$ is close to that of $I_{j'}$.
One concludes that
there exists a uniform $\gamma$ such that $0\,<\,\gamma\,<\,1$
$$|S(j_1..j_n)\,-\,S(j'_1..j'_m)|\,\leq\,\gamma^{\sharp(j\,\cap\,j')}
\;\;(inequality\,1)\;\;(2)$$
Here $\sharp(j\,\cap\,j')$ denotes the length of $j\cap\,j'$.
Therefore
for any infinite sequence $j\,=\,(j_1j_2...)$ the scaling function
$S$:
$$S(j)\,=\,\lim_{n\,\rightarrow\,\infty}\,S(j_1..j_n)$$
is well defined and has a H{\" o}lder modulus of continuity:
$$|S(j)\,-\,S(j')|\,\leq\,\gamma^{\sharp(j\cap\,j')}$$
This scaling function is canonically defined on the dual Cantor set
$\Sigma_d^{dual}$, whose elements are infinite sequences $(j_1j_2..)$.
Each such sequences should be thought of as a prescribed sequence of
inverse
branches of the shiftmap.\par
Say that a map is $C^{1+}$ if it is $C^{1\,+\,\epsilon}$ for some
$\epsilon$. \vskip .2 in
{\bf Theorem:} [Sullivan] {\it Every $C^{1+}$ embedding has a
H{\" o}lder continuous scaling
function. The scaling function is a $C^{1}$ invariant. Every
H{\" o}lder continuous function on the dual Cantor set with values in
$Simp_{2d-2}$ is the scaling function of a $C^{1+}$ embedding.}
\vskip .2 in
Here the H{\"o}lder continuity of the scaling function is defined with
respect to a metric on $\Sigma_d^{dual}$:
$$\rho_{\delta}(j,\,j')\,=\,exp(-\delta\,\sharp(j\cap j')$$
In the theorem $\delta$ (the metric on $\Sigma_d^{dual}$)
is not specified so we cannot specify $\epsilon$.

The problem which remained was to understand which functions occur
as scaling functions for $C^{1+1}$ and higher smoothness. Here we give
necessary and sufficient conditions for a function $S$ to arise
as a scaling function for a $C^{k\,+\,\epsilon}$ ($k$ positive integer
and $0\,<\,\epsilon\,\leq\,1$) embedding. The main observation is
that given an embedding, we should be able to extend
the identification map between $I_j$ and $I_{j'}$
to their convex hulls $<I_j>$ and $<I_{j'}$ to be $C^{k+\epsilon}$
close to affine provided $j\,\cap\,j'$ is long. Here close to affine is
measured after affinely rescaling
$<I_j>$ and $<I_{j'}>$ to the unit interval.
We refer to the process of changing the map by rescaling domain and
range to the
unit interval as renormalization.\vskip .2 in
We will first characterize those functions which are scaling functions
of $C^{1+\epsilon}$ Cantor sets. This is a special case of the main
theorem.
We state it seperately because of its simpler form.
Given a function $S:\,\Sigma_d^{dual}\,\rightarrow\,Simp_{2d-2}$.
We replace an arbitrary metric $\rho_{\delta}$ on $\Sigma_d^{dual}$
with a metric $\rho_S$ so that for an embedding with S as scaling
function there
exits $K$ so that for every $j,\,j'$:
$${{1}\over{K}}\,\leq\,{{|I_{j\,\cap\,j'}|}\over{\rho_S(j,\,j')}}
\,\leq\,K\;\;(3)$$
This metric is defined as:
$$\rho_S(j,\,j')\,=\,\sup_{w}\Pi_{t=1}^{n=\sharp(j\cap\,j')}\,
S(j_{t+1}j_{t+2}..j_nw)_{j_t}
$$
$(3)$ holds by $(1)$ because any infinite tail $w$ changes the product
by a uniformly bounded factor (by $(2)$).
\vskip .2 in
{\bf Theorem 1:} {\it Fix $0\,<\,\epsilon\,\leq\,1$.
The following are equivalent:\par
{\bf 1.} There exists a $C^{1+\epsilon}$ embedding with scaling
function
$S$.\par
{\bf 2.} S is $C^{\epsilon}$ on $(\Sigma_d^{dual},\,\rho_S)$. (Here
$C^1$ means
Lipschitz).}\par
{\bf Proof:} That {\bf 1.} $\Rightarrow$ {\bf 2.} follows when one
observes that a stronger form of $(2)$ holds:
$$|S(j_1..j_n)\,-\,S(j'_1..j'_n)|\,\leq\,K\,|I_{j\,\cap\,j'}|^{\epsilon}\;\;(4)
$$
This inequality carries over to the scaling function. Next apply
$(3)$.\par
That {\bf 2.} $\Rightarrow$ {\bf 1.}, i.e. the construction of a
$C^{1+\epsilon}$ Cantor set will be done in the
proof of the Main Theorem.

\QED
\vskip .2 in


{\bf Example 1:} For every $0<\e_1<\e_2\le 1$ there exists $S$
admitting a $C^{1+\e_1}$ embedding but not  $C^{1+\e_2}$. We find it as
follows:
Fon an arbitrary $0<\nu<{{\e_2-\e_1}\over 2}$ we can easily find a
function $S$ to
$Simp_{2d-2}$ which is $C^{1+\e_1+\nu}$ but not $C^{1+\e_2-\nu}$
on $\Sigma^{dual}_{d,n}$ with a standard metric $\rho_\delta$, $\delta
>\log d$.
We can find in fact $S$ so that for every $j\in \Sigma^{dual}_{d,n}$,
$i=1,..,d$
$|-\log S(j)_i/\delta\ -1|<\nu/\e_2$. This is chosen so that $S$ is
$C^{\e_1}$ but not $C^{\e_2}$ with respect to the metric $\rho_S$.

\

\

We now turn to the more intricate case of higher smoothness.\par
Let $A_1$ and $A_2$ be two subsets of the unit
interval $I\,=\,[0,\,1]$ such that both sets contain the
endpoints of $I$ and both have equal cardinality.  Denote the $k-th$
derivative operator by $D^k$ and denote by $D^k(A_1,\,A_2)$
the space of $C^k$ diffeomorphisms on $I$ which map $A_1$ to $A_2$.
For every constant $M\,>\,0$ consider the
space of $C^k$-diffeomorphisms:
$$D^k_{var}(M)(A_1,\,A_2)\,=\,\lbrace
\,\phi\,\in\,D^k(A_1,A_2):
\sup\,|D^k\phi(x)\,-\,D^k\phi(y)|\,<\,M
\rbrace$$
\vskip .2 in
{\bf Lemma:} {\it Assume that $A_1$ and $A_2$ consist of $2d$ points.
Assume
that $k\,<\,2\,d$. Then for each $f$, $g$ in $D^k_{var}(M)(A_1,A_2)$ we
have for
all integers $t\,\leq\,k$:
$$sup\,|D^tf\,-\,D^tg|\,\leq\,M$$}\par
{\bf Proof:} Consider two such maps $f$ and $g$. Their difference
vanishes on $A_1$. Since $2\,d\,>\,k$,
there exists (mean value theorem) for
each $t$ a point $x_t$ in $I$ for which:
$$D^tf(x_t)\,-\,D^tg(x_t)\,=\,0$$
The lemma follows by induction and integration.
\QED
\vskip .2 in
Given a function $S$ as above and a point $j$ in $\Sigma_d^{dual}$.
Consider
$S(j)$. It encodes a partition of $I$ in $ 2 d \,-\,1$ intervals.
Denote by
$A(j)$ the $2\,d$ end points of these intervals.
Consider any $j_0\,=\,1,..,d$ and consider the point
$j_0j$ in $\Sigma_d^{dual}$. Then $S(j_0j)$ specifies how the $j_0-th$
interval
in $j$ is subdivided. Consider two points $j$ and $j'$ in
$\Sigma_d^{dual}$
Every element in $D^k(A(j),A(j'))$ maps the $j_0-th$ interval in the
domain to the $j_0-th$ interval in the range, which we again can
renormalize.
This defines a map (restrict to $j_0-th$ interval and renormalize):
$$R_{j_0}:\,D^k(A(j),\,A(j'))\,\rightarrow\,D^k(\lbrace0,\,1\rbrace,\,
\lbrace0,\,1\rbrace)$$
\vskip .2 in
{\bf Main Theorem:} {\it Suppose $k\,<\,2\,d$. Suppose that we are
given a function $S$ as above. The following are
equivalent.\par
{\bf 1.} There exists a $C^{k+\epsilon}$ embedding with scaling
function
S.\par
{\bf 2.} There exists a constant $C$ so that for all $j$ and
$j'$ in $\Sigma_d^{dual}$ and all $j_0\,=\,1,..,d$:
$$D^k_{var}(C(j_0j,\,j_0j'))(A(j_0j),A(j_0j'))\,\cap\,
R_{j_0}(D^k_{var}(C(j,\,j'))(A(j),\,A(j'))\,\neq\,\emptyset$$
where for all $j,\,j'\,\in\,\Sigma_d^{dual}$,
$$C(j,\,j')\,=\,C\,\rho_S(j,\,j')^{k+\epsilon-1}$$}\vskip .2 in
{\bf Discussion of statement of theorem:}
The statement of the theorem may appear obscure. We briefly discuss
in an informal manner
how the scaling function records smoothness
beyond $C^1$. \par

{\bf 1)} Consider two strings $j$ and $j'$ and the identification
map between $I_j$ and $I_{j'}$. The scalings $S(j)$ and $S(j')$ record
how $2\,d$ specific points in $I_j$ map to $2\,d$ specific
points in $I_{j'}$. Consider the renormalized identification map
, and assume that we know that the variation of the
$k^{th}$ derivative of this identification map is small.
Consider any
$k\,+\,1$ of the $2\,d\,$ specific points. Since we know where
these points map, we can compute a value of the $k^{th}$ derivative
(just as the standard mean value theorem computes a value of the
first derivative given 2 points and their values).
Because the variation of the
$k-th$ derivative is small
, we obtain combinatorial relations
between any two choices of $k\,+\,1$ points.
Condition {\bf 2.} of the theorem  captures
this idea. It omits attempts to describe the derivative
algebraically\par

{\bf 2)} In fact we do not need all 2d points which appear in the
definition of the ratio geometry to be involved in the definition of
$D^k_{var}$'s, k+1 would be enough (see Lemma). In particular for
$C^{1+\e}$ the condition 2. makes impression we do not need the
geometry at all. However then the condition (4)in Prof of Theorem 1 is
hidden in 2. . Without (4) a map $I\to I$ in
$D^k_{var}(A(j_0j),A(j_0j'))$, even linear, after renormalizing  by
$R_{j_0}^{-1}$ may happen not to be extendible to a map belonging to
the second $D$ in 2. .\par

{\bf 3)} The condition of the main theorem seems to imply
that high smoothness is not discussed when $d$ is small.
We can however replace $d$ by any positive power $d^n$ in the following
manner. $\Sigma_d$ is canonically homeomorphic to $\Sigma_{d^n}$,
by the homeomorphism which groups the
digits of a point in $\Sigma_d$ in groups of $n$ digits. This
homeomorphism conjugates the $n-th$ iterate of the shiftmap on
$\Sigma_d$ to the shiftmap on $\Sigma_{d^n}$.\vskip .2 in
{\bf Proof of Main Theorem:}
We first show that {\bf 1} implies {\bf 2.}. Assume
that we are given a $C^{k\,+\epsilon}$ embedding $h$. Denote the
induced
shiftmap on the image by $f$. We may assume that its $d$
right-inverses extend as $C^{k\,+\,\epsilon}$ contractions to
the unit interval, the convex hull of the image of $h$. 
Denote by $f_{j'|j}$ the identification
between $<I_j>$ and $<I_{j'}>$ and denote by $F_{j'|j}$ the renormalized
identification defined on the unit interval $J$.
Then $f_{j'|j}$, respectively
$F_{j'|j}$, factors as a composition:
$$f_{j'|j}\,=\,f_{j'|j'\cap j}\,\circ\,f_{j'\cap j|j}$$
$$F_{j'|j}\,=\,F_{j'|j'\cap j}\,\circ\,F_{j'\cap j|j}$$
Since $f_{j'|j'\cap j}: <I_{j'\cap j}>\,\rightarrow\,<I_{j'}>$ is a
composition of $C^{k+\epsilon}$
contractions the derivatives of $f_{j'|j'\cap j}$ are controlled by the
first derivative.

More precisely, by a standard computation which we leave to the
reader,
there exists a constant $C$ so that for all $j$ and
$j'$,  all $1\,\leq\,t\,\leq\,k\,+\,\epsilon$
$$|f_{j'|j'\cap j}|_{t}\,\leq\,C\,|f_{j'|j'\cap j}|_1$$
Here $|.|_t$ denotes the supnorm of the $t-th$
derivative for $t$ integer and the $\alpha$- H{\" o}lder
norm of the $n-th$ derivative if $t\,=\,n\,+\,\alpha$,
$0\,<\,\alpha\,\,\leq\,1$.

But then:
$$|F_{j'|j'\cap j}|_t\,=\,{{|I_{j'\cap j}|^t}\over{|I_{j'}|}}\,
|f_{j'|
j'\cap j}|_t$$
$$\,\leq\,|I_{j'\cap j}|^{t-1}\,C$$
The last inequality follows because:
$${{|I_{j'}|}\over{|I_{j'\cap j}|}}\,=\,Df_{j'|j'\cap j}(x)$$
for some point $x\,\in\,I_{j'\cap j}$ and the bounded nonlinearity of
the maps.\par
Now let $j$ and $j'$ be two distinct points in $\Sigma_d^{dual}$.
Denote
by $j_n$, respectively $j_n'$ the beginning strings of length $n$.
Then for $n$ large enough $j\cap\,j'\,=\,j_n\,\cap\,j_n'$
and the sequence of maps $\lbrace\,F_{j'_n|j'\cap\,j}\rbrace$
is $C^{k\,+\,\epsilon}$-equicontinuous.
Since moreover:
$$F_{j'_{n+m}|j'\cap j}\,=\,F_{j'_{n+m}|j'_n}\,\circ\,F_{j'_n|j'\cap
j}$$
this sequence of maps is in fact $C^{k+\epsilon}$ convergent. Denote by
$F_{j'|j'\cap\,j}$
the limit map. By the same argument $F_{j|j'\cap\,j}$ is
defined. Therefore:
the limiting map:
$$F_{j'|j}\,=\,F_{j'|j'\cap\,j}\,\circ\,F_{j|j'\cap\,j}^{-1}$$
is well-defined and $C^{k+\epsilon}$ and therefore in $D_{var}^k$.
Since $\rho_S(j,\,j')$ is uniformly comparable
to $|I_{j'\cap\,j}|$ we obtain that this limiting map
$F_{j'|j}$ in $D_{var}^k(C')$ for some uniform constant $C'$. Since
moreover:
$$R_{j_0}\,F_{j'|j}\,=\,F_{j_0j'|j_0j}$$
we automatically have an element in the intersection.
 {\bf 2.} now follows.\vskip .2 in
We next show that {\bf 2.} implies {\bf 1.}. Since $S$ is given,
we first construct an embedding of the Cantor set with $S$ as scaling
function. We then show that this embedding is $C^{k\,+\,\epsilon}$.\par
Fix an arbitrary infinite word $w$. Construct a Cantor set $C$ in
the unit interval $<I>$ by consecutively subdividing any interval
$<I_j>$ according to $S(jw)$. We obtain an embedding with scaling
function
$S$. Denote the induced shiftmap on the image by $f_0$. It
is defined on a Cantor set $C$. In order
to show that this shiftmap has a $C^{k+\epsilon}$ extension, we
verify the assumptions to Whitney's extension theorem [Stein]. We will
construct functions $f_1,....f_k$ on $C$ so that for all
$x$, $y$ in $C$ and $l\,=\,0,...,k$ (Whitney conditions):
$$f_l(y)\,=\,\sum_{t=l}^{t=k}\,{{1}\over{(t-l)!}}\,f_{t-l}(x)(y-x)^{t-l}\,+\,
O(|y-x|^{k-l+\epsilon})$$
These functions $f_1,..f_k$ play the role of the first $k$ derivatives
of $f_0$.\par
The interval $<I>$ is subdivided in $d$ intervals $<I_{i}>$,
$i\,=\,1,..,d$. On each of the
intervals, $f_0$ maps $I_i\,=\,C\,\cap\,<I_i>$ to $I\,=\,
C\,\cap\,J$ by $f_0$. Now fix a
$i\,=\,1,..d$. We will work on each $<I_{i}>$ separately. For each
$t\,=\,1,..,k$ define $f_t$ on $I_{i}$ as:
$$f_t\,=\,lim_{n\,\rightarrow\,\infty}\,\lbrace\,D^t\phi_{j_1..j_n,j}\,
\rbrace_{(j_1..j_n)}$$
Here $\phi_{j_1..j_n,i}:\,J_{j_1..j_ni}\,\rightarrow\,J_{j_1..j_n}$
is any map whose renormalization is in
$$D_k^{var}(C(j_1..j_niw,\,j_1..j_nw)
(A(j_1..j_niw),\,A(j_1..j_nw))$$
We need to see that $f_t$ is in fact well-defined on the Cantor set.
We first verify that $f_t$ is defined point wise on the Cantor set.
Consider a string $j_1..j_n$ and an element $j_0$. For
$x\,\in\,I_{j_0j_1..j_ni}$, consider
$\phi_{j_0j_1..j_ni}(x)$ and $\phi_{j_1..j_ni}(x)$ and their
$t-th$ derivatives.
Then by assumption {\bf 2.} and the Lemma:
$$|D^t\phi_{j_0j_1..j_ni}(x)\,-\,D^t\phi_{j_1..j_ni}(x)|
\,\leq\,$$
$$
{{|I_{j_0j_1..j_ni}|}\over{|I_{j_0j_1..j_n}|^t}}\,C\,\rho_S
(j_1..j_niw,j_1..j_nw)^{k+\epsilon-1}\,\leq\,$$
$$C\,|I(j_1..j_ni)|^{k+\epsilon-t}$$
Therefore we obtain the convergence on the Cantor set
in fact exponentially fast.
\par
We need to check that the Whitney conditions hold on the
Cantor set. Let $x$ and $y$ be distinct points in the Cantor set in
$<I_i>$. Consider the first time that they wind up in different
intervals
in the subdivision:
$$x\,\in\,J_{j_0j_1..j_ni},\;\;y\,\in\,J_{j_0'j_1....j_ni},\;\;j_0\,\neq\,j_0'$$
Then again by {\bf 2.}:
$$|\phi_{j_1..j_ni}(y)\,-\,\phi_{j_1..j_ni}(x)\,-\,
\sum_{t=0}^{t=k}\,{{1}\over{t!}}\,D^t\phi_{j_1..j_ni}(x)(y-x)^t|\,\leq\,C\,
|x-y|^{k+\epsilon}$$
(and similar for the higher derivatives) where $C$ is a uniform
constant.
Since
$|D^t\phi_{j_1..j_ni}(x)\,-\,f_t(x)|\,\leq\,C\,|x\,-\,y|^{k+\epsilon-t}$,

we can take limits and obtain the Whitney conditions for the family
$f_0,f_1,...,f_k$. Consequently there exists a $C^{k+\epsilon}$
extension of $f$ to each $J_i$ and we have produced a $C^{k+\epsilon}$
embedding of $\Sigma_d$ with scaling function $S$.
\QED

\vskip .2 in
We say that two  embeddings
$h_1$ and $h_2$ are $C^r$ equivalent if the composition
$h_2\,\circ\,h_1^{-1}$ admits an extension as a $C^1$-diffeomorphism
to $\R$. It is well known that if $h_1$ and $h_2$ are $C^r$ equivalent
then the composition $h_2\,\circ\,h_1^{-1}$ in fact
admits an extension as a $C^r$-diffeomorphism. This result can also
be deduced as a corollary of the method employed in the main theorem.
\vskip .2 in

{\bf Corollary:} {\it Assume that $h_1$ and $h_2$ are equivalent
$C^{k\,+\,\epsilon}$ embeddings: $h_2\circ\,h_1^{-1}$ is $C^1$.
Then $h_2\circ\,h_1^{-1}$ is $C^{k+\e}$. }\par

{\bf Proof:}  To show that the conjugacy $h_2\circ\,h_1^{-1}$
has a $C^{k+\epsilon}$ extension,
it suffices to construct its higher derivatives on the Cantor set
and apply the Whitney extension theorem. This can be achieved
using the same manner as that employed in the second half of the proof
of the main theorem. Both embeddings have the same scaling function $S$
so , as the embeddings are $C^{k+1}$ the ratio geometries on finite
levels are close to one another in the sense of condition 2. of Main
Theorem.
\QED
{\bf Remark:} The preceding theorem is not totally satisfactory,
since we do not  understand how to extract $C^k$-smoothness
($k$ integer!)
from the scaling function. This is because in the previous scheme
everything which needs to be controlled is dominated by geometric
series. More refined finite smoothness categories
like $C^{1+zygmund}$ can however be treated in much the same way.
\vskip .2 in
We finally show in an example that conditions {\bf 2.} of the main
theorem can be explicitly checked, by constructing for
every $k,\,d,\,\epsilon$ with $k\,<\,2\,d\,-\,1$ and
$0\,<\,\epsilon\,<\,1$ an example of a scaling function with
a $C^{k\,+\,\epsilon}$ realization and none of higher degree of
smoothness.\par
{\bf Example:} Let $J_i\,=\,[{{2\,i\,-\,2}\over{2\,d\,-\,1}},
{{2\,i\,-\,1}\over{2\,d\,-\,1}}],\,i=1,..,d$

Define $f\,:
\cup_i\,J_i\,\rightarrow\,J\,=\,[0,1]$
as:
$$f(x)\,=\,A\,((2\,d\,-\,1)\,x\,+\,x^{k+\epsilon}),\;\;x\,\in\,J_1$$
while $f$ is affine on each $J_i,\;i\,\geq\,2$. Here the the constant
$A$ is chosen so that $f(J_1)\,=\,J$.

Of course the resulting Cantor set is $C^{k+\epsilon}$. We will
show that its scalingfunction on the dual Cantor set
has no $C^{k\,+\,\epsilon_1}$ realization
for all $\epsilon_1\,>\,\epsilon$, by explicitly checking
that condition {\bf 2.} of the main theorem does not hold for $k\,+\,
\epsilon_1$.

Let $w$ be any element
in $\Sigma_d^{dual}$ which does not contain the symbol $1$. Denote
by $1_n$ the string of length $n$ consisting of $1's$ only:
$$1_n\,=\,11..1$$
Consider the infinite strings $j\,=\,1_nw$ and
$j'\,=\,1_n1w\,=\,1_{n+1}w$.
Consider the subdivision $A(j)$, respectively $A(j')$, of the unit
interval
dictated by $S(j)$ and $S(j')$.
Let $\Phi_n$ be any map in
$D_{var}^k(A(j),\,A(j')$ for which its renormalized restriction
$R_1\Phi$ is in fact in $D^k_{var}(A(1j),\,A(1j')$.
We will bound the variation of the $k-th$
derivative of $\Phi_n$ from below and
conclude that condition {\bf 2.} of the main theorem is not satisfied
 with $k\,+\,\epsilon_1$.

We denote by $A(j)_m$ the $m-th$ point
from the left in $A(j)$.
Because $k\,<\,2\,d\,-\,1$
there exists $x\,\in\,[A(j)_2,\,A(j)_{2d}]$ such that:
$$D^kF_{j'|j}(x)\,=\,D^k\,\Phi_n(x)$$
Recall that $F_{j'|j}$ is the renormalization of $f_{j'|j}$ for the map
$f$ defined above. See the notation of the proof of the Main Theorem.
Similarly there exists $y\,\in\,[A(1j)_1,..,A(1j)_{2d-1}]$ so that:
$$D^kF_{1j'|1j}(y)\,=\,D^k(R_1\Phi_n)(y)$$.
We have that:
$$D^kF_{j'|j}(x)\,=\,B\,\,{2d-1}^{-n(k-1+\e)}\,x^{\epsilon}$$
(note that $(2d-1)^{-n}\,\sim\,\rho_S(j',j)$).
The map $F_{1j'|1j}$ is just the renormalization of the
restriction of the limit map $F_{j'|j}$ to the left most interval
$[A(j)_1,A(j)_2]$ in the unit interval. Let $y'$ be the point in the
interval $[A(j)_1,A(j)_2]$, corresponding to $y$ after rescaling the
unit
interval back to $[A(j)_1,A(j)_2]$. Then we have that:
$$D^k(F_{j'|j})(y')\,=\,B\,\,(2d-1)^{-n(k-1+\e)}\,(y')^{\epsilon}$$
where $B$ is a computable constant.

But $|x\,-\,y'|\,>\,const\,(2d-1)^2$. Consequently:
$$D^k\Phi_n(x)\,-\,D^k\Phi_n(y')\,=\,const\,\,(2d-1)^{-n(k-1+\e)}\,(
x^{\epsilon}\,-\,(y')^{\epsilon}$$
and is comparable to:
$$\rho_s(j',j)^{k-1\,+\,\epsilon}$$
i.e. the variation of $D^k\Phi_n$ is at least on the order of:
$\rho_s(j',j)^{k-1\,+\,\epsilon}$.

Since:
$$\lim_{n\,\rightarrow\,\infty}\,
{{\rho_s(j',j)^{k-1\,+\,\epsilon_1}}\over{\rho_s(j',j)^{k-1\,+\,\epsilon}}}
\,=\,0$$
condition {\bf 2.} of the theorem can not be satisfied for
$$C(j',\,j)\,=\,C\,\rho_S(j',j)^{k-1\,+\,\epsilon_1}$$

\vskip .3 in
{\bf Reference:}\par
[Stein] Singular Integrals\par

[Sullivan] Weyl proceedings AMS.

\vfill\eject\end